\documentclass[letterpaper, 10pt]{article}
\usepackage{amsthm,amsmath,amsfonts}

\newtheorem{thm}{Theorem}[section]
\newtheorem{cor}[thm]{Corollary}

\newtheorem{prop}[thm]{Proposition}

\theoremstyle{definition}

\newtheorem{rem}[thm]{Remark}
\newtheorem{conj}[thm]{Conjecture}

\newcommand{\mR}{\mathcal{R}}

\title{A note on regular Ramsey graphs}
\author{
Noga Alon\thanks{School of Computer Science and School of
Mathematical Sciences, Raymond and Beverly Sackler Faculty of Exact
Sciences, Tel Aviv University, Tel Aviv 69978, Israel.  E-mail:
nogaa@post.tau.ac.il. Research supported in part by an ERC advanced
grant, by the Israel Science Foundation and by a USA-Israel BSF
grant.} \and Sonny Ben-Shimon\thanks{School of Computer Science,
Raymond and Beverly Sackler Faculty of Exact Sciences, Tel Aviv
University, Tel Aviv 69978, Israel. E-mail: sonny@post.tau.ac.il.
Research conducted as part of the author's Ph.D. thesis under the
supervision of Prof. Michael Krivelevich.} \and Michael
Krivelevich\thanks{School of Mathematical Sciences, Raymond and
Beverly Sackler Faculty of Exact Sciences, Tel Aviv University, Tel
Aviv 69978, Israel. E-mail: krivelev@post.tau.ac.il. Research
supported in part by USA-Israel BSF Grant 2006322, by grant 1063/08
from the Israel Science Foundation, and by a Pazy memorial award.} }
\begin{document}
\maketitle
\begin{abstract}
We prove that there is an absolute constant $C>0$ so that
for every natural $n$ there exists a
triangle-free \emph{regular} graph with no independent set of size
at least $C\sqrt{n\log n}$.
\end{abstract}
\section{Introduction}
A major problem in extremal combinatorics asks to determine
the maximal $n$ for which there exists a graph $G$ on $n$ vertices
such that $G$ contains no triangles and no independent set of size
$t$. This Ramsey-type problem was settled asymptotically by Kim
\cite{Kim95} in 1995, after a long line of research; Kim showed that
$n=\Theta(t^2/\log t)$. Recently,
Bohman \cite{Boh2009} gave an alternative proof of Kim's result by
analyzing the so-called triangle-free process, as proposed by
Erd\H{o}s, Suen and Winkler \cite{ErdSueWin95}, which is a natural
way of generating a triangle-free graph. Consider now the above
problem with the additional constraint that $G$ must be regular. In
this short note we show that the same asymptotic results hold up to
constant factors. The main ingredient of the proof is a gadget-like
construction that transforms a triangle-free graph with no
independent set of size $t$, which is not too far from being regular,
into a triangle-free \emph{regular} graph with no independent set of
size $2t$.

Our main result can be stated as follows.
\begin{thm}\label{t:main}
There is a positive constant $C$ so that
for every natural $n$
there exists a regular triangle-free graph $G$ on $n$ vertices whose
independence number satisfies
$\alpha(G) \leq C \sqrt{n\log n}$.
\end{thm}
Denote by $R(k,\ell)$ the maximal $n$ for which there exists a graph
on $n$ vertices which contains neither a complete subgraph
on $k$ vertices nor an
independent set on $\ell$ vertices. Let $R^{\mathrm{reg}}(k,\ell)$
denote the maximal $n$ for which there exists a \emph{regular} graph
on $n$ vertices which contains neither a complete subgraph
on $k$ vertices nor
an independent set on $\ell$ vertices. Clearly, for every $k$ and
$\ell$ one has $R^{\mathrm{reg}}(k,\ell)\leq R(k,\ell)$.
Theorem \ref{t:main} states that
$R^{\mathrm{reg}}(3,t)=\Theta\left(R(3,t)\right)
=\Theta\left(\frac{t^2}{\log
t}\right)$.
\section{Proof of Theorem \ref{t:main}}
Note first that
the statement of the theorem is trivial for small values of
$n$. Indeed, for every $n_0$
one can choose the constant $C$ in the theorem so that for
$ n \leq n_0$, $C \sqrt {n \log n} \geq n$, implying that for such values of
$n$ a graph with no edges satisfies the assertion of the theorem.
We thus
may and will assume, whenever this is needed during the proof,
that $n$ is sufficiently large.

The following well known theorem due to Gale and to Ryser gives a
necessary and sufficient condition for two lists of non-negative
integers to be the degree sequences of the classes of vertices
of a simple bipartite graph. The proof
follows easily from the max-flow-min cut condition on the
appropriate network flow graph (see e.g. \cite[Theorem
4.3.18]{West2001}).
\begin{thm}[Gale; Ryser]\label{t:GayRys57}
If $\mathbf{d}=(d_1,\ldots,d_m)$ and
$\mathbf{d'}=(d'_1,\ldots,d'_n)$ are lists of non-negative integers
with $d_1\geq\ldots\geq d_m$, $d'_1\geq\ldots\geq d'_n$ and $\sum
d_i= \sum d'_j$ then there exists a simple bipartite graph with
degree sequences $\mathbf{d}$ and $\mathbf{d'}$ on each side
respectively iff $\sum_{i=1}^m \min\{d_i,s\}\geq \sum_{j=1}^s d'_j$
for every $1\leq s\leq n$ .
\end{thm}

\begin{cor}\label{c:noga1}
Let $a \geq 1$ be a real.
If $\mathbf{d}=(d_1,\ldots,d_m)$ is a list of non-negative integers
with $d_1\geq\ldots\geq d_m$ and
\begin{equation}
\label{e21}
d_1\leq\min\left\{ad_m,
\frac{4am}{(a+1)^2}\right\},
\end{equation}
then there exists a simple bipartite graph with degree sequence
$\mathbf{d}$ on each side. In particular, this holds for
$d_1 \leq \min \{2d_m, \frac{8m}{9}\}. $
\end{cor}
\begin{proof}
By Theorem  \ref{t:GayRys57} it suffices to check that for every
$s$, $ 1 \leq s \leq m$, $\sum_{i=1}^s d_i \leq
\sum_{i=1}^m \min\{d_i,s\}$. Suppose this is not the case
and there is some $s$ as above so that
\begin{equation}
\label{e22}
d_1+ d_2 + \ldots +d_s > \sum_{i=1}^m \min\{d_i,s\}.
\end{equation}
If $d_i<d_1$ for some $i$ satisfying $2 \leq i \leq s$, replace
$d_i$ by $d_1$. Observe that by doing so the left hand side of
(\ref{e22}) increases by $d_1-d_i$, whereas the right hand side
increases by at most this quantity, hence (\ref{e22}) still holds
with this new value of $d_i$. We can thus assume that
$d_1=d_2 = \cdots =d_s$. Note that if $d_1 \leq s$, then
(\ref{e22}) cannot hold, hence $d_1 >s$. If $d_i >d_1/a$ for some $i$
satisfying $s+1 \leq i \leq m$, then reducing
it to $d_1/a$ (even if this is not an integer), maintains
(\ref{e22}), as the left hand side does not change,
whereas the right hand side  can only decrease. Moreover, the new
sequence still satisfies (\ref{e21}). Thus we may assume that
in (\ref{e22})
$d_i =d_1/a$ for all $s+1 \leq i \leq m$.
Put $d=d_{i+1}~(= d_{i+2}= \ldots =d_m)$, then (\ref{e22})
gives
$$
d_1 + \ldots + d_s =s\cdot(ad) > \sum_{i=1}^m \min\{d_i,s\}
=s^2+(m-s)d.
$$
Therefore $[(a+1)s-m]d >s^2$, implying that
$(a+1)s-m>0$, that is, $s > \frac{m}{a+1}$, and
\begin{equation}
\label{e23}
d> \frac{s^2}{(a+1)s-m}.
\end{equation}
The function $g(s) =\frac{s^2}{(a+1)s-m}$ attains its minimum
in the range $ \frac{m}{a+1} < s \leq m$ at $s=\frac{2m}{a+1}$
and its value at this point is $\frac{4m}{(a+1)^2}.$
We thus conclude from (\ref{e23}) that
$d > \frac{4m}{(a+1)^2}$ and hence that
$d_1=ad > \frac{4am}{(a+1)^2}$ contradicting the assumption
(\ref{e21}). This completes the proof.
\end{proof}
\begin{rem}
Condition (\ref{e21}) is tight for all values of $a>1$, in the sense
that if $d_1> \frac{4am}{(a+1)^2}$ and $d_1=d_2 \ldots =d_s$
for $s=\frac{2m}{a+1}$ with $d_i=\frac{4m}{(a+1)^2}$
for all $s+1 \leq i \leq m$, then there is no simple bipartite graph
whose degree sequence in  each  side is $(d_1, d_2, \ldots ,d_m)$.
This follows from Theorem \ref{t:GayRys57}.
\end{rem}

Let $\mR(n,3,t)$ denote the set of all triangle-free graphs $G$ on
$n$ vertices with $\alpha(G)<t$. As usual, let $\Delta(G)$ and
$\delta(G)$ denote the respective maximal and minimal degrees of $G$.
\begin{prop}\label{p:noga1}
Let $t$ and $d$ be integers. If there exists a graph
$G\in\mR(n,3,t)$ such that $\Delta(G)-\delta(G)\leq d\leq
\frac{4}{9}\cdot\left\lfloor\frac{n}{\Delta(G)+1}\right\rfloor$,
then there exists a $(d+\Delta(G))$-\emph{regular} graph
$G'\in\mR(2n,3,2t-1)$.
\end{prop}
\begin{proof}
Construct a new graph $G'$ as follows. Take two copies of $G$, and
color each of these copies by the same equitable coloring using
$\Delta(G)+1$ colors with all color classes of cardinality either
$\left\lfloor n/(\Delta(G)+1)\right\rfloor$ or $\left\lceil
n/(\Delta(G)+1)\right\rceil$ using the Hajnal-Szemer\'{e}di Theorem
\cite{HajSze70} (see also a shorter proof due to Kierstead and
Kostochka \cite{KieKos2008}). Let $C$ and $C'$ be the same color
class in each of the copies of $G$. Denote the degree sequence of
the vertices of $C$ in $G$ by $d'_1\leq\ldots\leq d'_m$, where
$m=|C|$, and set $d_i=d+\Delta(G)-d'_i$. According to Corollary
\ref{c:noga1} there exists a simple bipartite graph with $m$
vertices on each side, where the degree sequence of each side is
$d_1\geq\ldots\geq d_m$ as the maximal degree
$d_1=d+\Delta(G)-\delta(G)\leq 2d$, the minimal degree $d_m\geq d$,
and by our assumption on $G$ we have $d_1\leq \frac{8m}{9}$. We can
thus connect the vertices of $C$ and $C'$ using this bipartite graph
such that all vertices in $C\cup C'$ have degree $d+\Delta(G)$. By
following this method for every color class, we create the graph
$G'$ which is $(d+\Delta(G))$-regular, triangle-free and has no
independent set of cardinality $2t-1$.
\end{proof}

\subsection{The $H$-free process and Bohman's result}
Consider the following randomized greedy algorithm to generate a
graph on $n$ labeled vertices with no $H$-subgraph for some fixed
graph $H$. Given a set of $n$ vertices, a sequence of graphs
$\{G^{(H)}_i\}^t_{i=0}$ on this set of vertices is constructed.
Start with $G^{(H)}_0$ as the empty graph, and for each $0<i\leq t$,
the graph $G^{(H)}_i$ is defined by $G^{(H)}_{i-1}\cup\{e_i\}$ where
$e_i$ is chosen uniformly at random from all unselected pairs of
vertices that do not create a copy of $H$ when added to
$G^{(H)}_{i-1}$. The process terminates at step $t$, the first time
that no potential unselected pair $e_{t+1}$ exists. This algorithm
is called the \emph{$H$-free process}.

The $K_3$-free process was proposed by Erd\H{o}s, Suen and Winkler
\cite{ErdSueWin95} and was further analyzed by Spencer \cite{Spe95}.
Recently, Bohman \cite{Boh2009} extending and improving previous
results, was able to analyze the $K_3$-free process and to show that
with high probability it passes through an almost regular
Ramsey-type graph.
\begin{thm}[Bohman \cite{Boh2009}]\label{t:Boh2009}
With high probability\footnote{In this context we mean that the
mentioned events hold with probability tending to $1$ as $n$, the
number of vertices, goes to infinity.} there exists an integer
$1\leq m=m(n)$ such that the following properties hold
simultaneously:
\begin{enumerate}
\item $G^{(K_3)}_m\in\mR(n,3,C\sqrt{n\log n})$ for some absolute constant $C>0$;
\item $\Delta(G^{(K_3)}_m)=\Theta(\sqrt{n\log n})$;
\item\label{i:degdiff} $\Delta(G^{(K_3)}_m)-\delta(G^{(K_3)}_m)=o(\sqrt{n/ \log n})$.
\end{enumerate}
\end{thm}
\begin{rem}
Item \eqref{i:degdiff} can be derived implicitly from
\cite{Boh2009}, or alternatively, it follows from \cite[Theorem
1.4]{BohKeePre}, as the degree of every vertex is a \emph{trackable
extension variable}.
\end{rem}

Note that Proposition \ref{p:noga1} in conjunction with Theorem
\ref{t:Boh2009} completes the proof of Theorem \ref{t:main} for
every large enough \emph{even} integer $n$. To fully complete the
proof, we describe how to deal with the case of $n$ odd. So, let now
$n$ be be large enough and odd. Our aim is to show the existence of
a regular triangle-free graph $G_n$ on $n$ vertices with
$\alpha(G_n)=O(\sqrt{n\log n})$. The approach we take to achieve
this goal is to construct a ``big'' graph satisfying our Ramsey
conditions on an even number of vertices, and to add to it a
``small'' graph with an odd number of vertices without affecting the
asymptotic results claimed.

For every $k=0\pmod{5}$, and every even $r\leq 2k/5$, let $H_{k,r}$
denote a graph constructed as follows. Start with a copy of $C_5$
blown up by factor of $k/5$ and delete from the resulting graph
$(2k/5 - r/2)$ disjoint $2$-factors (which exist by Petersen's
Theorem, see e.g. \cite[Theorem 3.3.9]{West2001}). $H_{k,r}$ is
hence a triangle-free $r$-regular graph on $k$ vertices.

Denote by $F_n$ an $r$-regular triangle-free graph on $2n$ vertices
with $\alpha(F_n)\leq C \sqrt{n\log n}$ for some absolute constant
$C$, and furthermore assume $r$ is even (this can be achieved by
choosing the appropriate parameter $d$ in Proposition \ref{p:noga1},
as we have much room to spare with the values we plug in from
Theorem \ref{t:Boh2009}). Let $n_0=(n-k)/2$, where $k=5\pmod{10}$,
and $k=(1+o(1))\frac{5C}{2}\sqrt{n\log n}$. Clearly, $n_0$ is
integer. The graph $F_{n_0}$ is $r$-regular for some even $r\leq
\alpha(F_{n_0})$, is triangle-free on $2n_0$ vertices, and satisfies
$\alpha(F_{n_0})\leq C\sqrt{n_0\log n_0}\leq C\sqrt{n\log n}$. Now,
define $G_n$ to be a disjoint union of $F_{n_0}$ and $H_{k,r}$.
Clearly, $G_n$ has $2n_0+k=n$ vertices, is $r$-regular,
triangle-free and satisfies
$\alpha(G_n)=\alpha(F_{n_0})+\alpha(H_{k,r})\le \alpha(F_{n_0})+k
\leq C\sqrt{n\log n} + k =O(\sqrt{n\log n})$.
\section{Discussion}
A natural question that extends the above is to try and determine
$R^{\mathrm{reg}}(k,\ell)$ for other values of $k$ and $\ell$ (in
particular for fixed values of $k>3$ and $\ell$ arbitrary large),
and also to try and investigate its relation with $R(k,\ell)$.
The following conjecture seems plausible.
\begin{conj}
\label{c31}
For every $k \geq 2$ there is a constant $c_k>0$ so that
$R^{\mathrm{reg}}(k,\ell) \geq c_k R(k,\ell)$ for all $\ell \geq 2$.
\end{conj}
This is trivial for $k=2$,  and  by our main result here
holds for $k=3$  as well.

Recently, Bohman and Keevash \cite{BohKeePre} were able to
generalize the techniques of \cite{Boh2009} for the $H$-free
process, where $H$ is a strictly 2-balanced graph. This in turn
provided new lower bounds for $R(k,\ell)$ (as complete graphs are
strictly 2-balanced) where $k$ is fixed and $\ell$ arbitrarily
large. It is plausible to think that these results can also be used
to construct \emph{regular} Ramsey graphs in a manner similar to
that described in this note. Nonetheless, since the asymptotic
behavior of $R(k,\ell)$ is not known for $k \geq 4$, a complete
proof of Conjecture \ref{c31} appears to require some additional
ideas, and remains open.

\bibliographystyle{abbrv}
\bibliography{regramsey}
\end{document}